\documentstyle[12pt]{article}
\begin{document}

\author{S.V. Ludkovsky.}

\title{Wrap groups of fiber bundles over quaternions and octonions.}

\date{15 January 2008}
\maketitle

\begin{abstract}
This article is devoted to the investigation of wrap groups of
connected fiber bundles over the fields of real $\bf R$, complex
$\bf C$ numbers, the quaternion skew field $\bf H$ and the octonion
algebra $\bf O$. These groups are constructed with mild conditions
on fibers. Their examples are given. It is shown, that these groups
exist and for differentiable fibers have the infinite dimensional
Lie groups structure, that is, they are continuous or differentiable
manifolds and the composition $(f,g)\mapsto f^{-1}g$ is continuous
or differentiable depending on a class of smoothness of groups.
Moreover, it is demonstrated that in the cases of real, complex,
quaternion and octonion manifolds these groups have structures of
real, complex, quaternion or octonion manifolds respectively.
Nevertheless, it is proved that these groups does not necessarily
satisfy the Campbell-Hausdorff formula even locally.
\end{abstract}

\section{Introduction.} Wrap groups of fiber bundles considered
in this paper are constructed with the help
of families of mappings from a fiber bundle with a marked point into
another fiber bundle with a marked point over the fields $\bf R$,
$\bf C$, $\bf H$ and the octonion algebra $\bf O$. Conditions on
fibers supplied with parallel transport structures are rather mild
here. Therefore, they generalize geometric loop groups of circles,
spheres and fibers with parallel transport structures over them. A
loop interpretation is lost in their generalizations, so they are
called here wrap groups. This paper continues previous works of the
author on this theme, where generalized loop groups of manifolds
over $\bf R$, $\bf C$ and $\bf H$ were investigated, but neither for
fibers nor over octonions \cite{ludan,lugmlg,lujmslg,lufoclg}.

\par Loop groups of circles were first introduced by Lefshetz in
1930-th and then their construction was reconsidered by Milnor in
1950-th. Lefshetz has used the $C^0$-uniformity on families of
continuous mappings, which led to the necessity of combining his
construction with the structure of a free group with the help of
words. Later on Milnor has used the Sobolev's $H^1$-uniformity, that
permitted to introduce group structure more naturally
\cite{milmorse}. Iterations of these constructions produce iterated
loop groups of spheres. Then their constructions were generalized
for fibers over circles and spheres with parallel transport
structures over $\bf R$ or $\bf C$ \cite{gaj}.

\par Wrap groups of quaternion and octonion fibers
as well as for wider classes of fibers over $\bf R$ or $\bf C$ are
defined and investigated here for the first time.

\par Holomorphic functions of quaternion and octonion variables were
investigated in \cite{luoyst,luoyst2,lufsqv}. There specific
definition of super-differentiability was considered, because the
quaternion skew field has the graded algebra structure. This
definition of super-differentiability does not impose the condition
of right or left super-linearity of a super-differential, since it
leads to narrow class of functions. There are some articles on
quaternion manifolds, but practically they undermine a complex
manifold with additional quaternion structure of its tangent space
(see, for example, \cite{museya,yano} and references therein).
Therefore, quaternion manifolds as they are defined below were not
considered earlier by others authors (see also \cite{lufsqv}).
Applications of quaternions in mathematics and physics can be found
in \cite{emch,guetze,hamilt,lawmich}.

\par In this article wrap groups of different classes of smoothness are
considered. Henceforth, we consider not only orientable manifolds
$M$ and $N$, but also nonorientable manifolds.
\par  In particular,
geometric loop groups have important applications in modern physical
theories (see \cite{ish,mensk} and references therein). Groups of
loops are also intensively used in gauge theory. Wrap groups defined
below with the help of families of mappings from a manifold $M$ into
another manifold $N$ with a dimension $dim (M) >1$ can be used in
the membrane theory which is the generalization of the string
(superstring) theory.
\par Section 2 is devoted to the definitions of topological and
manifold structures of wrap groups. The existence of these groups is
proved and that they are infinite dimensional Lie groups not
satisfying even locally the Campbell-Hausdorff formula (see Theorems
3, 6, 12, Corollaries 5, 8, 9 and Examples 10). In the cases of
complex, quaternion and octonion manifolds it is proved that they
have structures of complex, quaternion and octonion manifolds
respectively.
\par All main results of this paper are obtained for the first time.

\section{Wrap groups of fibers.}
To avoid misunderstandings we first give our definitions and
notations.

\par {\bf 1.1. Note.} Denote by ${\cal A}_r$ the
Cayley-Dickson algebra such that ${\cal A}_0=\bf R$, ${\cal A}_1=\bf
C$, ${\cal A}_2=\bf H$ is the quaternion skew field, ${\cal A}_3=\bf
O$ is the octonion algebra. Henceforth we consider only $0\le r\le
3$.

\par {\bf 1.2. Definition.}
~ A canonical closed subset $Q$ of the Euclidean space $X=\bf R^n$
or of the standard separable Hilbert space $X=l_2({\bf R})$ over
$\bf R$ is called a quadrant if it can be given by the condition
$Q:=\{ x\in X: q_j(x)\ge 0 \} $, where $(q_j: j\in \Lambda _Q)$ are
linearly independent elements of the topologically adjoint space
$X^*$. Here $\Lambda _Q\subset \bf N$ (with $card (\Lambda _Q)=k\le
n$ when $X=\bf R^n$) and $k$ is called the index of $Q$. If $x\in Q$
and exactly $j$ of the $q_i$'s satisfy $q_i(x)=0$ then $x$ is called
a corner of index $j$.
\par If $X$ is an additive group and also left and right module over
$\bf H$ or $\bf O$ with the corresponding associativity or
alternativity respectively and distributivity laws then it is called
the vector space over $\bf H$ or $\bf O$ correspondingly.

In particular $l_2 ({\cal A}_r)$ consisting of all sequences $x = \{
x_n\in {\cal A}_r: n \in {\bf N} \} $ with the finite norm $\| x \|
<\infty $ and scalar product $(x,y):=\sum_{n=1}^{\infty } x_ny_n^*$
with $\| x \| := (x,x)^{1/2}$ is called the Hilbert space (of
separable type) over ${\cal A}_r$, where $z^*$ denotes the
conjugated Cayley-Dickson number, $zz^* =: |z|^2$, $z\in {\cal
A}_r$. Since the unitary space $X={\cal A}_r^n$ or the separable
Hilbert space $l_2({\cal A}_r)$ over ${\cal A}_r$ while considered
over the field $\bf R$ (real shadow) is isomorphic with $X_{\bf
R}:=\bf R^{2^rn}$ or $l_2({\bf R})$, then the above definition also
describes quadrants in ${\cal A}_r^n$ and $l_2({\cal A}_r)$. In the
latter case we also consider generalized quadrants as canonical
closed subsets which can be given by $Q:=\{ x\in X_{\bf R}:$
$q_j(x+a_j)\ge 0, a_j\in X_{\bf R}, j\in \Lambda _Q \} ,$ where
$\Lambda _Q\subset \bf N$ ($card(\Lambda _Q)=k\in \bf N$ when
$dim_{\bf R}X_{\bf R}<\infty $).

\par {\bf 1.2.2. Definition.}
A differentiable mapping $f:  U\to U'$ is called a diffeomorphism if
\par $(i)$ $f$ is bijective and there exist continuous mappings
$f'$ and $(f^{-1})'$, where $U$ and $U'$ are interiors of quadrants
$Q$ and $Q'$ in $X$.
\par In the ${\cal A}_r$ case with $1\le r\le 3$ we consider
bounded generalized quadrants $Q$ and $Q'$ in ${\cal A}_r^n$ or
$l_2({\cal A}_r)$ such that they are domains with piecewise
$C^{\infty }$-boundaries. We impose additional conditions on the
diffeomorphism $f$ in the $1\le r\le 3$ case:
\par $(ii)$ ${\bar \partial }f=0$ on $U$,
\par $(iii)$ $f$ and all its strong (Frech\'et) differentials
(as multi-linear operators) are bounded on $U$, where $\partial f$
and ${\bar \partial }f$ are differential $(1,0)$ and $(0,1)$ forms
respectively, $d=\partial +{\bar \partial }$ is an exterior
derivative, for $2\le r\le 3$ $\partial $ corresponds to
super-differentiation by $z$ and ${\tilde \partial }={\bar \partial
}$ corresponds to super-differentiation by ${\tilde z}:=z^*$, $z\in
U$ (see \cite{luoyst,luoyst2}).

\par The Cauchy-Riemann Condition $(ii)$ means that $f$ on $U$ is the
${\cal A}_r$-holomorphic mapping.

\par {\bf 1.2.3. Definition and notation.}
An ${\cal A}_r$-manifold $M$ with corners is defined in the usual
way:  it is a metric separable space modelled on $X={\cal A}_r^n$ or
$X=l_2({\cal A}_r)$ respectively and is supposed to be of class
$C^{\infty }$, $0\le r\le 3$. Charts on $M$ are denoted $(U_l, u_l,
Q_l)$, that is, $u_l:  U_l\to u_l(U_l) \subset Q_l$ is a $C^{\infty
}$-diffeomorphism for each $l$, $U_l$ is open in $M$, $u_l\circ
{u_j}^{-1}$ is biholomorphic for $1\le r\le 3$ from the domain
$u_j(U_l\cap U_j)\ne \emptyset $ onto $u_l(U_l\cap U_j)$ (that is,
$u_j\circ u_l^{-1}$ and $u_l\circ u_j^{-1}$ are holomorphic and
bijective) and $u_l\circ u_j^{-1}$ satisfy conditions $(i-iii)$ from
\S 1.2.2, $\bigcup_jU_j=M$.

\par A point $x\in M$ is called a corner of index $j$
if there exists a chart $(U,u,Q)$ of $M$ with $x\in U$ and $u(x)$ is
of index $ind_M(x)=j$ in $u(U)\subset Q$. A set of all corners of
index $j\ge 1$ is called a border $\partial M$ of $M$, $x$ is called
an inner point of $M$ if $ind_M(x)=0$, so $\partial M=\bigcup_{j\ge
1}\partial ^jM$, where $\partial ^jM:=\{ x\in M:  ind_M(x)=j \} $.
\par For a real manifold with corners on the connecting mappings
$u_l\circ u_j^{-1}\in C^{\infty }$ of real charts only
Condition $1.2.2(i)$ is imposed.

\par {\bf 1.2.4. Terminology.}
In an ${\cal A}_r$-manifold $N$ there exists an
Hermitian metric, which in each analytic system of coordinates is
the following $\sum_{j,k=1}^nh_{j,k}dz_jd{\bar z}_k$, where
$(h_{j,k})$ is a positive definite Hermitian matrix with
coefficients of the class $C^{\infty }$, $h_{j,k}=h_{j,k}(z)\in
{\cal A}_r$, $z$ are local coordinates in $N$.
\par As real manifolds we shall consider Riemann manifolds.

\par In accordance with the definition above for internal points of
$N$ it is supposed that they can belong only to interiors of charts,
but for boundary points $\partial N$ it may happen that $x\in
\partial N$ belongs to boundaries of several charts. It is
convenient to choose an atlas such that $ind (x)$ is the same for
all charts containing this $x$.

\par {\bf 1.3.1. Remark.} If $M$ is a metrizable space and $K=K_M$
is a closed subset in $M$ of codimension $codim_{\bf R} ~N\ge 2$ such that
$M\setminus K =M_1$ is a manifold with corners over ${\cal
A}_r$, then we call $M$ a pseudo-manifold over ${\cal A}_r$,
where $K_M$ is a critical subset.
\par Two pseudo-manifolds $B$ and $C$ are called diffeomorphic, if
$B\setminus K_B$ is diffeomorphic with $C\setminus K_C$ as for
manifolds with corners (see also \cite{gaj,michor}).

\par Take on $M$ a Borel $\sigma $-additive measure $\nu $
such that $\nu $ on $M\setminus K$ coincides with the Riemann volume
element and $\nu (K)=0$, since the real shadow of $M_1$ has it.
\par The uniform space $H^t_p(M_1,N)$ of all continuous piecewise
$H^t$ Sobolev mappings from $M_1$ into $N$ is introduced in the
standard way \cite{lujmslg,lufoclg}, which induces $H^t_p(M,N)$ the
uniform space of continuous piecewise $H^t$ Sobolev mappings on $M$,
since $\nu (K)=0$, where ${\bf R}\ni t\ge [m/2]+1$, $m$ denotes the
dimension of $M$ over $\bf R$, $[k]$ denotes the integer part of
$k\in \bf R$, $[k]\le k$. Then put $H^{\infty }_p(M,N)=\bigcap_{t>m}
H^t_p(M,N)$ with the corresponding uniformity.
\par For manifolds over ${\cal A}_r$ with $1\le r\le 3$ take
as $H^t_p(M,N)$ the completion of the family of all continuous
piecewise ${\cal A}_r$-holomorphic mappings from $M$ into $N$
relative to the $H^t_p$ uniformity, where $[m/2]+1\le t\le \infty $.
Henceforth we consider pseudo-manifolds with connecting mappings of charts
continuous in $M$ and $H^{t'}_p$ in $M\setminus K_M$
for $0\le r\le 3$, where $t'\ge t$.
\par {\bf 1.3.2. Note.} Since the octonion algebra $\bf O$
is non-associative, we consider a non-associative subgroup $G$ of
the family $Mat_q({\bf O})$ of all square $q\times q$ matrices with
entries in $\bf O$. More generally $G$ is a group which has a
$H^t_p$ manifold  structure over ${\cal A}_r$ and group's operations
are $H^t_p$ mappings. The $G$ may be non-associative for $r=3$, but
$G$ is supposed to be alternative, that is, $(aa)b=a(ab)$ and
$a(a^{-1}b)=b$ for each $a, b\in G$.

\par As a generalization of pseudo-manifolds there is used the
following (over $\bf R$ and $\bf C$ see \cite{gaj,souriau}). Suppose
that $M$ is a Hausdorff topological space of covering dimension $dim
~M=m$ supplied with a family $ \{ h : U\to M \} $ of the so called
plots $h$ which are continuous maps satisfying conditions $(D1-D4)$:
\par $(D1)$ each plot has as a domain a convex subset $U$ in ${\cal
A}_r^n$, $n\in \bf N$;
\par $(D2)$ if $h: U\to M$ is a plot, $V$ is a convex subset in
${\cal A}_r^l$ and $g: V\to U$ is an $H^t_p$ mapping, then $h\circ
g$ is also a plot, where $t\ge [m/2]+1$;
\par $(D3)$ every constant map from a convex set $U$ in ${\cal
A}_r^n$ into $M$ is a plot;
\par $(D4)$ if $U$ is a convex set in ${\cal A}_r^n$ and
$ \{ U_j: j\in J \} $ is a covering of $U$ by convex sets in ${\cal
A}_r^n$, each $U_j$ is open in $U$, $h: U\to M$ is such that each
its restriction $h|_{U_j}$ is a plot, then $h$ is a plot. Then $M$ is
called an $H^t_p$-differentiable space.
\par A mapping $f: M\to N$ between two $H^t_p$-differentiable spaces
is called differentiable if it continuous and for each plot $h: U\to
M$ the composition $f\circ h: U\to N$ is a plot of $N$. A
topological group $G$ is called an $H^t_p$-differentiable group if
its group operations are $H^t_p$-differentiable mappings.
\par  Let $E$, $N$, $F$ be $H^{t'}_p$-pseudo-manifolds or
$H^{t'}_p$-differentiable spaces over ${\cal A}_r$, let also $G$ be
an $H^{t'}_p$ group over ${\cal A}_r$, $t\le t'\le \infty $. A fiber
bundle $E(N,F,G,\pi ,\Psi )$ with a fiber space $E$, a base space
$N$, a typical fiber $F$ and a structural group $G$ over ${\cal
A}_r$, a projection $\pi : E\to N$ and an atlas $\Psi $ is defined
in the standard way \cite{gaj,michor,sulwint} with the condition,
that transition functions are of $H^{t'}_p$ class such that for
$r=3$ a structure group may be non-associative, but alternative.

\par Local trivializations $\phi _j\circ \pi \circ \Psi _k^{-1}:
V_k(E)\to V_j(N)$ induce the $H^{t'}_p$-uniformity in the family $W$ of all principal $H^{t'}_p$-fiber bundles $E(N,G,\pi ,\Psi )$, where
$V_k(E) = \Psi _k(U_k(E))\subset X^2(G)$, $V_j(N) = \phi _j(U_j(N))\subset X(N)$, where $X(G)$ and $X(N)$ are ${\cal A}_r$-vector spaces on which $G$ and $N$ are modelled, $(U_k(E),\Psi _k)$ and $(U_j(N),\phi _j)$ are charts of
atlases of $E$ and $N$, $\Psi _k =\Psi _k^{E}$, $\phi _j = \phi _j^{N}$.

\par If $G=F$ and $G$ acts on itself by left shifts, then a fiber bundle is
called the principal fiber bundle and is denoted by $E(N,G,\pi ,\Psi )$. As a
particular case there may be $G={\cal A}_r^*$, where ${\cal A}_r^*$
denotes the multiplicative group ${\cal A}_r\setminus \{ 0 \} $.
If $G=F= \{ e \} $, then $E$ reduces to $N$.

\par {\bf 2. Definitions.} Let $M$ be a connected
$H^t_p$-pseudo-manifold over ${\cal A}_r$, $0\le r\le 3$ satisfying
the following conditions:

\par $(i)$ it is compact;

\par $(ii)$ $M$ is a union of two closed subsets
over ${\cal A}_r$ $A_1$ and $A_2$, which are pseudo-manifolds and which are canonical closed subsets in $M$ with $A_1\cap A_2=\partial A_1\cap
\partial A_2=:A_3$ and a codimension over $\bf R$ of $A_3$ in $M$ is
$codim_{\bf R}A_3=1$, also $A_3$ is a
pseudo-manifold;

\par $(iii)$ a finite set of marked points $s_{0,1},...,s_{0,k}$
is in $\partial A_1\cap \partial A_2$, moreover, $\partial A_j$ are
arcwise connected $j=1, 2$;

\par $(iv)$ $A_1\setminus \partial A_1$ and $A_2\setminus \partial A_2$
are $H^t_p$-diffeomorphic with $M\setminus [ \{ s_{0,1},...,s_{0,k}
\} \cup (A_3\setminus Int (\partial A_1\cap \partial A_2))]$ by
mappings $F_j(z)$, where $j=1$ or $j=2$, $\infty \ge t\ge [m/2]+1$,
$m=dim_{\bf R}M$ such that $H^t\subset C^0$ due to the Sobolev
embedding theorem \cite{miha}, where the interior $Int (\partial A_1\cap \partial A_2)$ is taken in $\partial A_1\cup \partial A_2$.

\par Instead of $(iv)$ we consider also the case
\par $(iv')$  $M$, $A_1$ and $A_2$ are such that
$(A_j\setminus \partial A_j)\cup \{ s_{0,1},...,s_{0,k} \} $ are \\
$C^0([0,1],H^t_p(A_j,A_j))$-retractable on $X_{0,q}\cap A_j$, where
$X_{0,q}$ is a closed arcwise connected subset in $M$, $j=1$ or
$j=2$, $s_{0,q}\in X_{0,q}$, $X_{0,q}\subset K_M$, $q=1,...,k$,
$codim_{\bf R}~ K_M\ge 2$.
\par Let $\hat M$ be a compact connected $H^t_p$-pseudo-manifold
which is a canonical
closed subset in ${\cal A}_r^l$ with a boundary $\partial {\hat M}$
and marked points $\{ {\hat s}_{0,q}\in \partial {\hat M}:
q=1,...,2k \} $ and an $H^t_p$-mapping $\Xi : {\hat M}\to M$ such
that \par $(v)$ $\Xi $ is surjective and bijective from ${\hat
M}\setminus \partial {\hat M}$ onto $M\setminus \Xi (\partial {\hat
M})$ open in $M$, $\Xi ({\hat s}_{0,q})=\Xi ({\hat s}_{0,k+q}) =
s_{0,q}$ for each $q=1,...,k$, also $\partial M\subset \Xi (\partial
{\hat M})$.

\par A parallel transport structure on a $H^{t'}_p$-differentiable
principal $G$-bundle $E(N,G,\pi ,\Psi )$ with arcwise connected $E$
and $G$ for $H^t_p$-pseudo-manifolds $M$ and $\hat M$ as above over
the same ${\cal A}_r$ with $t'\ge t+1$ assigns to each $H^t_p$
mapping $\gamma $ from $M$ into $N$ and points $u_1,...,u_k\in
E_{y_0}$, where $y_0$ is a marked point in $N$, $y_0=\gamma
(s_{0,q})$, $q=1,...,k$, a unique $H^t_p$ mapping ${\bf P}_{{\hat
\gamma },u}: {\hat M}\to E$ satisfying conditions $(P1-P5)$: \par
$(P1)$ take ${\hat \gamma }: {\hat M}\to N$ such that ${\hat \gamma
}=\gamma \circ \Xi $, then ${\bf P}_{{\hat \gamma },u}({\hat s}_{0,q})=u_q$
for each $q=1,...,k$ and $\pi \circ {\bf P}_{{\hat \gamma },u}={\hat
\gamma }$
\par $(P2)$ ${\bf P}_{{\hat \gamma },u}$ is the $H^t_p$-mapping
by $\gamma $ and $u$;
\par $(P3)$ for each $x\in \hat M$ and every $\phi \in DifH^t_p
({\hat M}, \{ {\hat s}_{0,1},...,{\hat s}_{0,2k} \} )$ there is the
equality ${\bf P}_{{\hat \gamma },u}(\phi (x)) = {\bf P}_{{\hat
\gamma }\circ \phi ,u}(x)$, where $DifH^t_p({\hat M}, \{ {\hat
s}_{0,1},...,{\hat s}_{0,2k} \} )$ denotes the group of all $H^t_p$
homeomorphisms of $\hat M$ preserving marked points $\phi ({\hat
s}_{0,q})={\hat s}_{0,q}$ for each $q=1,...,2k$;
\par $(P4)$ ${\bf P}_{{\hat \gamma },u}$ is $G$-equivariant,
which means that ${\bf P}_{{\hat \gamma },uz}(x) = {\bf P}_{{\hat
\gamma },u}(x)z$ for every $x\in {\hat M}$ and each $z\in G$;
\par $(P5)$ if $U$ is an open neighborhood of ${\hat s}_{0,q}$ in
$\hat M$ and ${\hat \gamma }_0, {\hat \gamma }_1: U\to N$ are
$H^{t'}_p$-mappings such that ${\hat \gamma }_0({\hat
s}_{0,q})={\hat \gamma }_1({\hat s}_{0,q})=v_q$ and tangent spaces,
which are vector manifolds over ${\cal A}_r$, for $\gamma _0$ and
$\gamma _1$ at $v_q$ are the same, then the tangent spaces of ${\bf
P}_{{\hat \gamma }_0,u}$ and ${\bf P}_{{\hat \gamma }_1,u}$ at $u_q$
are the same, where $q=1,...,k$, $u=(u_1,...,u_k)$.

\par Two $H^{t'}_p$-differentiable principal $G$-bundles
$E_1$ and $E_2$ with parallel transport structures $(E_1,{\bf P}_1)$
and $(E_2,{\bf P}_2)$ are called isomorphic, if there exists an
isomorphism $h: E_1\to E_2$ such that ${\bf P}_{2,{\hat \gamma
},u}(x) = h ({\bf P}_{1,{\hat \gamma }, h^{-1}(u)}(x))$ for each
$H^t_p$-mapping $\gamma : M\to N$ and $u_q\in (E_2)_{y_0}$, where
$q=1,...,k$, $h^{-1}(u)=(h^{-1}(u_1),...,h^{-1}(u_k))$.

\par Let $(S^ME)_{t,H}:=
(S^{M, \{ s_{0,q}: q=1,...,k \} }E; N,G,{\bf P})_{t,H}$ be a set of
$H^t_p$-closures of isomorphism classes of $H^t_p$ principal $G$
fiber bundles with parallel transport structure.

\par {\bf 3. Theorems.} {\it  {\bf 1.} The uniform space
$(S^ME)_{t,H}$ from \S 2 has the structure of a topological
alternative monoid with a unit and with a cancelation property and
the multiplication operation of $H^l_p$ class with $l=t'-t$
($l=\infty $ for $t'=\infty $). If $N$ and $G$ are separable, then
$(S^ME)_{t,H}$ is separable. If $N$ and $G$ are complete, then
$(S^ME)_{t,H}$ is complete. \par {\bf 2.} If $G$ is associative,
then $(S^ME)_{t,H}$ is associative. If $G$ is commutative, then
$(S^ME)_{t,H}$ is commutative. If $G$ is a Lie group, then
$(S^ME)_{t,H}$ is a Lie monoid. \par {\bf 3.} The $(S^ME)_{t,H}$ is
non-discrete, locally connected and infinite dimensional for
$dim_{\bf R}(N\times G)>1$.}

\par {\bf Proof.} If there is a homomorphism $\theta : G\to F$
of $H^{t'}_p$-differentiable groups, then there exists an induced
principal $F$ fiber bundle $(E\times ^{\theta }F)(N,F,\pi ^{\theta
}, \Psi ^{\theta })$ with the total space $(E\times ^{\theta }F)=
(E\times F)/{\cal Y}$, where $\cal Y$ is the equivalence relation
such that $(vg,f){\cal Y} (v,\theta (g)f)$ for each $v\in E$, $g\in
G$, $f\in F$. Then the projection $\pi ^{\theta }: (E\times ^{\theta
}F)\to N$ is defined by $\pi ^{\theta }([v,f]) =\pi (v)$, where
$[v,f] := \{ (w,b): (w,b){\cal Y} (v,f), w\in E, b\in F \} $ denotes
the equivalence class of $(v,f)$. \par Therefore, each parallel
transport structure $\bf P$ on the principal $G$ fiber bundle
$E(N,G,\pi ,\Psi )$ induces a parallel transport structure ${\bf
P}^{\theta }$ on the induced bundle by the formula ${\bf P}^{\theta
}_{{\hat \gamma }, [u,f]}(x) = [{\bf P}_{{\hat \gamma },u}(x),f]$.

\par Define multiplication with the help of certain embeddings and
isomorphisms of spaces of functions. Mention that for each two
compact canonical closed subsets $A$ and $B$ in ${\cal A}_r^l$
Hilbert spaces $H^t(A,{\bf R}^m)$ and $H^t(B,{\bf R}^m)$ are
linearly topologically isomorphic, where $l, m\in \bf N$, hence
$H^t_p(A,N)$ and $H^t_p(B,N)$ are isomorphic as uniform spaces. Let
$H^t_p(M,\{ s_{0,1},...,s_{0,k} \} ;W,y_0) := \{ (E,f): E=E(N,G,\pi
,\Psi )\in W, f={\bf P}_{{\hat \gamma },y_0}\in H^t_p: \pi \circ
f(s_{0,q})=y_0 \forall q=1,...,k; \pi \circ f={\hat \gamma }, \gamma
\in H^t_p(M,N) \} $ be the space of all $H^{t'}_p$ principal $G$
fiber bundles $E$ with their parallel transport $H^t_p$-mappings
$f={\bf P}_{{\hat \gamma },y_0}$, where $W$ is as in \S 1.3.2.
Put $\omega _0=(E_0,{\bf P}_0)$
be its element such that $\gamma _0(M)= \{ y_0 \} $, where $e\in G$
denotes the unit element, $E_0=N\times G$, $\pi _0(y,g)=y$ for each
$y\in N$, $g\in G$, ${\bf P}_{{\hat \gamma }_0,u}={\bf P}_0$.
\par The mapping $\Xi : {\hat M}\to M$ from \S 2 induces the embedding
\par $\Xi ^*: H^t_p(M,\{ s_{0,1},...,s_{0,k} \} ;W,y_0)\hookrightarrow
H^t_p({\hat M},\{ {\hat s}_{0,1},...,{\hat s}_{0,2k} \} ;W,y_0)$, \\
where $\hat M$ and ${\hat A}_1$ and ${\hat A}_2$ are retractable
into points.

\par Let as usually $A\vee B := \rho ({\cal Z})$ be the wedge sum of
pointed spaces $(A,\{ a_{0,q}: q=1,...,k \} )$ and $(B,\{ b_{0,q}:
q=1,...,k \} )$, where ${\cal Z} := [A\times \{ b_{0,q}: q=1,...,k
\}\cup \{ a_{0,q}: q=1,...,k \}\times B]\subset A\times B$, $\rho $
is a continuous quotient mapping such that $\rho (x)=x$ for each
$x\in {\cal Z}\setminus \{ a_{0,q}\times b_{0,j}; q, j=1,...,k \} $
and $\rho (a_{0,q})=\rho (b_{0,q})$ for each $q=1,...,k$, where $A$
and $B$ are topological spaces with marked points $a_{0,q}\in A$ and
$b_{0,q}\in B$, $q=1,...,k$. Then the wedge product $g\vee f$ of two
elements $f, g\in H^t_p(M, \{ s_{0,1},...,s_{0,k} \} ;N,y_0)$ is
defined on the domain $M\vee M$ such that $(f\vee g)(x\times
b_{0,q})=f(x)$ and $(f\vee g)(a_{0,q}\times x)=g(x)$ for each $x\in
M$, where to $f, g$ there correspond $f_1, g_1 \in H^t_p({\hat M},
\{ {\hat s}_{0,1},...,{\hat s}_{0,2k} \} ;N,y_0)$ such that $f_1 =
f\circ \Xi $ and $g_1 = g\circ \Xi $.

\par Let $(E_j,{\bf P}_{{\hat \gamma }_j,u^j})\in
H^t_p(M,\{ s_{0,1},...,s_{0,k} \} ;W,y_0)$, $j=1, 2$, then take
their wedge product ${\bf P}_{{\hat \gamma },u^1} := {\bf P}_{{\hat
\gamma }_1,u^1}\vee {\bf P}_{{\hat \gamma }_2,v}$ on $M\vee M$ with
$v_q=u_q g_{2,q}^{-1}g_{1,q+k}= y_0\times g_{1,q+k}$ for each
$q=1,...,k$ due to the alternativity of $G$, $\gamma =\gamma _1\vee
\gamma _2$, where ${\bf P}_{{\hat \gamma }_j,u^j}({\hat s}_{j,0,q})=
y_0\times g_{j,q}\in E_{y_0}$ for every $j$ and $q$. For each
$\gamma _j: M\to N$ there exists ${\tilde \gamma }_j: M\to E_j$ such
that $\pi \circ {\tilde \gamma }_j=\gamma _j$. Denote by ${\bf m}: G\times
G\to G$ the multiplication operation. The wedge product $(E_1,{\bf
P}_{{\hat \gamma }_1,u^1})\vee (E_2,{\bf P}_{{\hat \gamma }_2,u^2})$
is the principal $G$ fiber bundle $(E_1\times E_2)\times ^{\bf m}G$
with the parallel transport structure ${\bf P}_{{\hat \gamma
}_1,u^1}\vee {\bf P}_{{\hat \gamma }_2,v}$.

\par The uniform space
$H^t_p(J,A_3;W,y_0):=\{ (E,f)\in H^t_p(J,W): \pi \circ f(A_3)=\{ y_0
 \} \} $ has the $H^t_p$-manifold structure and has an
embedding into \\ $H^t_p(M,\{ s_{0,1},...,s_{0,k} \} ;W,y_0)$ due to
Conditions 2$(i-iii)$, where either $J=A_1$ or $J=A_2$. This induces
the following embedding $\chi ^*: H^t_p(M\vee M,\{ s_{0,q}\times
s_{0,q}: q=1,...,k \} ;W,y_0)\hookrightarrow H^t_p(M,\{ s_{0,q}:
q=1,...,k \} ;W,y_0)$. \par Analogously considering $H^t_p(M,\{
X_{0,q}: q=1,...,k \} ;W,y_0)=\{ f\in H^t(M,W): f(X_{0,q})=\{ y_0
\}, q=1,...,k \} $ and $H^t_p(J,A_3\cup \{ X_{0,q}: q=1,...,k \}
;W,y_0)$ in the case $(iv')$ instead of $(iv)$ we get the embedding
$\chi ^*: H^t_p(M\vee M,\{ X_{0,q}\times X_{0,q}: q=1,...,k \}
;W,y_0)\hookrightarrow H^t_p(M,\{ X_{0,q}: q=1,...,k \} ;W,y_0)$.
Therefore, $g\circ f:=\chi ^*(f\vee g)$ is the composition in
$H^t_p(M,\{ s_{0,q}: q=1,...,k \} ;W,y_0)$.
\par There
exists the following equivalence relation $R_{t,H}$ in $H^t_p(M,\{
X_{0,q}: q=1,...,k \} ;W,y_0)$: $fR_{t,H}h$ if and only if there
exist nets $\eta _n\in DifH^t_p(M, \{ X_{0,q}: q=1,..., k \} )$,
also $f_n$ and $h_n\in H^t_p(M,\{ X_{0,q}: q=1,...,k \} ;W,y_0)$
with $\lim_n f_n=f$ and $\lim_n h_n=h$ such that $f_n(x)=h_n(\eta
_n(x))$ for each $x\in M$ and $n\in \omega $, where $\omega $ is a
directed set and convergence is considered in $H^t_p(M,\{ X_{0,q}:
q=1,...,k \} ;W,y_0)$. Henceforward in the case 2$(iv)$ we get
$s_{0,q}$ instead of $X_{0,q}$ in the case 2$(iv')$.

\par Thus there exists the quotient uniform space \\
$H^t_p(M,\{ X_{0,q}: q=1,...,k \} ; W,y_0)/R_{t,H} =: (S^ME)_{t,H}$.
In view of \cite{omo,omori2} $DifH^t_p(M)$ is the group of
diffeomorphisms for $t\ge [m/2]+1$. The Lebesgue measure $\lambda $
in the real shadow of ${\hat M}$ by the mapping $\Xi $ induces the
measure $\lambda ^{\Xi }$ on $M$ which is equivalent to $\nu $,
since $\Xi $ is the $H^t_p$-mapping from the compact space onto the
compact space, $\lambda (\partial {\hat M})=0$ and $\Xi : {\hat
M}\setminus {\partial {\hat M}}\to M$ is bijective.
\par Due to Conditions $(P1-P5)$ each element $f={\bf P}_{{\hat \gamma
},u}$ up to a set $Q_M$ of measure zero, $\nu (Q_M)=0$, is given as
$f\circ \Xi ^{-1}$ on $M\setminus Q_M$, where $\pi \circ f={\hat
\gamma }$, ${\hat \gamma }=\gamma \circ \Xi $. Denote $f\circ \Xi
^{-1}$ also by $f$. Thus, for each $(E,f)\in H^t_p(M,\{ s_{0,q}:
q=1,...,k \} ;W,y_0)$ the image $f(M)$ is compact and connected in
$E$. \par Therefore, for each partition $Z$ there exists $\delta
>0$ such that for each partition $Z^*$ with $\sup_i \inf_j
dist(M_i,{M^* }_j) <\delta $ and $(E,f)\in H^t(M,W;Z)$,
$f(s_{0,q})=u_q$, there exists $(E,f_1)\in H^t(M,W;Z^*)$ with
$f_1(s_{0,q})=u_q$ for each $q=1,...,k$ such that $fR_{t,H}f_1$,
where $M_i$ and $M^*_j$ are canonical closed pseudo-submanifolds in
$M$ corresponding to partitions $Z$ and $Z^*$, $H^t(M,W;Z)$ denotes
the space of all continuous piecewise $H^t$-mappings from $M$ into
$W$ subordinated to the partition $Z$ such that $Z$ and $Z^*$
respect $H^t_p$ structure of $M$.

\par Hence there exists a countable subfamily
$\{ Z_j: j \in {\bf N} \} $ in the family of all partitions
$\Upsilon $ such that $Z_j\subset Z_{j+1}$ for each $j$ and $\lim_j
\tilde diam Z_j=0$. Then
\par $(i)$ $str-ind \{ H^t (M,\{ s_{0,q}: q=1,...,k \} ;W,y_0;Z_j);
h^{Z_i}_{Z_j}; {\bf N} \} /R_{t,H}=(S^ME)_{t,H}$ is separable if $N$ and $G$ are separable, since each space $H^t_p(M,\{ s_{0,q}: q=1,...,k \} ;W,y_0;Z_j)$ is separable.

\par The space $str-ind \{ H^t(M,\{ s_{0,q}: q=1,...,k \} ;W,y_0;Z_j);
h^{Z_i}_{Z_j}; {\bf N} \}$ is complete due to Theorem 12.1.4
\cite{nari}, when $N$ and $G$ are complete. Each class of
$R_{t,H}$-equivalent elements is closed in it. Then to each Cauchy
net in $(S^ME)_{t,H}$ there corresponds a Cauchy net in $str-ind \{
H^t(M\times [0,1],\{ s_{0,q}\times e\times 0;W,y_0;Z_j\times Y_j);
h^{Z_i \times Y_i}_{Z_j\times Y_j}; {\bf N} \}$ due to theorems
about extensions of functions \cite{miha,seel,touger}, where $Y_j$
are partitions of $[0,1]$ with $\lim_j\tilde diam(Y_j)=0$,
$Z_j\times Y_j$ are the corresponding partitions of $M\times [0,1]$.
Hence $(S^ME)_{t,H}$ is complete, if $N$ and $G$ are complete.
\par If $f, g\in H^t(M,X)$ and $f(M)\ne g(M)$, then
\par $(ii)$ $\inf_{\psi \in DifH^t_p(M, \{ s_{0,q}: q=1,...,k \})} \| f\circ
\psi - g \| _{H^t(M,X)}>0$. Thus equivalence classes $<f>_{t,H}$ and
$<g>_{t,H}$ are different.  The pseudo-manifold $\hat M$ is arcwise
connected. Take $\eta : [0,1]\to \hat M$ an $H^t_p$-mapping with
$\eta (0)={\hat s}_{0,q}$ and $\eta (1)={\hat s}_{0,k+q}$, where
$1\le q\le k$. Choose in $\hat M$ $H^t_p$-coordinates one of which
is a parameter along $\eta $. Therefore, for each $g_q, g_{k+q}\in
G$ there exists ${\bf P}_{{\hat \gamma },u}$ with ${\bf P}_{{\hat
\gamma },u}(s_{0,q})=y_0\times g_q$ and ${\bf P}_{{\hat \gamma
},u}(s_{0,k+q})=y_0\times g_{k+q}$ for each $q=1,...,k$. Since $E$
and $G$ are arcwise connected, then $N$ is arcwise connected and
$(S^ME)_{t,H}$ is locally connected for $dim_{\bf R}N>1$. Thus, the
uniform space $(S^ME)_{t,H}$ is non-discrete. \par The tangent
bundle $TH^t_p(M,E)$ is isomorphic with $H^t_p(M,TE)$, where $TE$ is
the $H^{t'-1}_p$ fiber bundle, $t'\ge t+1$. There is an infinite
family of $f_{\alpha }\in H^t_p(M,TE)$ with pairwise distinct images
in $TE$ for different $\alpha $ such that $f_{\alpha }(M)$ is not
contained in $\bigcup_{\beta <\alpha }f_{\beta }(M)$, $\alpha \in
\Lambda $, where $\Lambda $ is an infinite ordinal. Therefore,
$T(S^ME)_{t,H}$ is an infinite dimensional fiber bundle due to
$(ii)$ and inevitably $(S^ME)_{t,H}$ is infinite dimensional.
\par Evidently, if $f\vee g=h\vee g$ or $g\vee f=g\vee h$ for
$\{ f, g, h \} \subset H^t_p(M,\{ s_{0,q}: q=1,...,k \} ;W,y_0)$,
then $f=h$. Thus $\chi ^*(f\vee g)=\chi ^*(h\vee g)$ or $\chi
^*(g\vee f)=\chi ^*(g\vee h)$ is equivalent to $f=h$ due to the
definition of $f\vee g$ and the definition of equal functions, since
$\chi ^*$ is the embedding. Using the equivalence relation $R_{t,H}$
gives $<f>_{t,H}\circ <g>_{t,H}= <h>_{t,H}\circ <g>_{t,H}$ or
$<g>_{t,H}\circ <f>_{t,H}= <g>_{t,H}\circ <h>_{t,H}$ is equivalent
to $<h>_{t,H}=<f>_{t,H}$. Therefore, $(S^ME)_{t,H}$ has the
cancelation property.
\par Since $G$ is alternative, then
$a_{2,q}[a_{2,q}^{-1}(a_{2,q+k}(a_{2,q}^{-1}a_{1,q+k}))]=
a_{2,q+k}(a_{2,q}^{-1}a_{1,q+k})$, hence ${\bf P}_1\vee ({\bf
P}_2\vee {\bf P}_2)= ({\bf P}_1\vee {\bf P}_2)\vee {\bf P}_2$; also
$a_{2,q}[a_{2,q}^{-1}(a_{1,q+k}(a_{1,q}^{-1}a_{1,q+k}))]=
a_{1,q+k}(a_{1,q}^{-1}a_{1,q+k})$, consequently, ${\bf P}_1\vee
({\bf P}_1\vee {\bf P}_2)= ({\bf P}_1\vee {\bf P}_1)\vee {\bf P}_2$
and inevitably for equivalence classes $(aa)b=a(ab)$ and
$b(aa)=(ba)a$ for each $a, b\in (S^ME)_{t,H}$. Thus $(S^ME)_{t,H}$
is alternative.

\par If $G$ is associative, then the parallel transport structure
gives $(f\vee g)\vee h=f\vee (g\vee h)$ on $M\vee M\vee M$ for each
$\{ f, g, h \} \subset H^t_p(M,\{ s_{0,q}: q=1,...,k;W,y_0)$.
Applying the embedding $\chi ^*$ and the equivalence relation
$R_{t,H}$ we get, that $(S^ME)_{t,H}$ is associative $<f>_{\xi
}\circ (<g>_{\xi }\circ <h>_{\xi })= (<f>_{\xi }\circ <g>_{\xi
})\circ <h>_{\xi }$.

\par In view of Conditions 2$(i-iv)$ there exists an
$H^t_p$-diffeomoprhism of $(A_1\setminus A_3)\vee (A_2\setminus
A_3)$ with $(A_2\setminus A_3)\vee (A_1\setminus A_3)$ as pseudo-manifolds
(see \S 1.3.1). For the
measure $\nu $ on $M$ naturally the equality $\nu (A_3)=0$ is
satisfied. If $M'$ - is the submanifold may be with corners or
pseudo-manifold, accomplishing the partition $Z=Z_f$ of the manifold
$M$, then the codimension $M'$ in $M$ is equal to one and $\nu
(M')=0$. For the point $s_{0,q}$ in $(M \setminus A_3)\cup \{
s_{0,q} \} $ there exists an open neighborhood $U$ having the
$H^t_p$-retraction $F: [0,1]\times U\to \{ s_{0,q} \}$. Hence it is
possible to take a sequence of diffeomorphisms $\psi _n \in
DifH^t_p(M, \{ s_{0,q}: q=1,..., k \} )$ such that $\lim_{n\to
\infty } diam (\psi _n(U))=0$.
\par Let $w_0$ be a mapping $w_0: M\to W$ such that $w_0(M)=\{ y_0
\times e \} $. Consider $w_0\vee (E,f)$ for some $(E,f)\in H^t_p
(M,\{ s_{0,q}: q=1,...,k \} ;W,y_0)$. If $(E,f)\in H^t_p(M,\{
s_{0,q}: q=1,..,k \}; W,y_0)$ with the natural positive $t\in \bf
N$, then $f$ is bounded relative to the uniformity of the uniform
space $H^t_p(M;E)$. If $U_n$ is a sequence of bounded open or
canonical closed subsets in $M$ such that $\lim_n diam (U_n)=0$,
then $\lim_{n\to \infty } \nu (V_n)=0$ for the sequence of $\nu
$-measurable subsets $V_n$ such that $V_n\subset U_n$. Therefore,
for each bounded sequence $\{ g_n: g_n\in H^t_p(M;E); n\in {\bf N}
\} $ there exists the limit $\lim_{n\to \infty } g_n|_{U_n} =0$
relative to the $H^t_p$ uniformity, where $U_n$ is subordinated to
the partition of $M$ into $H^t$ submanifolds. Then if $\{ g_n:
g_n\in H^t_p(M,\{ s_{0,q}: q=1,...,k \} ;E,y_0); n\in {\bf N} \} $
is a bounded sequence such that $g_n$ converges to $g\in H^t_p(M,\{
s_{0,q}: q=1,...,k \} ;N,y_0)$ on $M\setminus W_k$ for each $k$
relative to the $H^t_p$-uniformity, the given open $W_k$ in $M$,
where $k, n\in \bf N$ and $\lim_{n\to \infty } \nu
(W_n\bigtriangleup U_n)=0$, then $g_n$ converges to $g$ in the
uniform space $H^t_p(M,\{ s_{0,q}: q=1,...,k \} ;E,y_0)$.

\par Mention that for each marked point $s_{0,q}$ in $M$ there exists
a neighborhood $U$ of $s_{0,q}$ in $M$ such that for each
$\gamma _1\in H^t_p(M, \{ s_{0,q}: q=1,...,k \}; N,y_0)$ there exists
$\gamma _2\in H^t_p$ such that they are $R_{t,H}$ equivalent and $\gamma _2|_U = y_0$. Therefore, if $C$ is an arcwise connected compact subset in $M$ of codimension $codim_{\bf R} C\ge 1$ such that $s_{0,q}\in C$, then the standard proceeding shows that for each $\gamma _1\in H^t_p$ there exists
$\gamma _2\in H^t_p$ such that $\gamma _1R_{t,H}\gamma _2$ and $\gamma _2|_C= y_0$. Since $C$ is compact, then each its open covering has a finite subcovering and hence \par $(Y_0)$ there exists an open neighborhood $U$ of $C$ in $M$ such that for each $\gamma _1$ there exists $\gamma _2$ such that
$\gamma _1R_{t,H}\gamma _2$ and $\gamma _2|_U = y_0$.

\par There
exists a sequence $\eta _n \in DifH^t_p(M, \{ s_{0,q}: q=1,...,k \}
)$ such that $\lim_{n\to \infty }diam(\eta _n(A_2\setminus \partial A_2))=0$ and $w_n, f_n\in H^t_p(M,\{ s_{0,q}: q=1,...,k \} ;E,y_0)$ with \par $(iii)$
$\lim_{n\to \infty }f_n=f$, $\lim_{n\to \infty }w_n=w_0$ and
$\lim_{n\to \infty }\chi ^*(f_n\vee w_n)(\eta _n^{-1})=f$ due to
$\pi \circ f(s_{0,q})=s_{0,q}$ in the formula of differentiation of
compositions of functions (over $\bf H$ and $\bf O$ see it in
\cite{luoyst,luoyst2,lufsqv}).
\par In more details, the sequence $\eta _n$ as a
limit of $\eta _n(A_2)$ produces a pseudo-submanifold $B$ in $M$ of codimension not
less than one such that $B$ can be presented with the help of the
wedge product of spheres and compact quadrants up to
$H^t_p$-diffeomorphism with marked points $\{ s_{0,q}: q=1,...,k \}
$, but as well $B$ may be a finite discrete set also. Then by induction the procedure can be continued lowering the
dimension of $B$. Particularly there may be circles and curves in
the case of the unit dimension. Two quadrants up to an
$H^t_p$ quotient mapping gluing boundaries produce a sphere. Thus the consideration reduces to
the case of the wedge product of spheres. The case of spheres
reduces to the iterated construction with circles, since the reduced
product $S^1\wedge S^n$ is $H^t_p$ homeomorphic with $S^{n+1}$ (see
Lemma 2.27 \cite{swit} and \cite{gaj}). For the particular case of
the $n$-dimensional sphere $M_n=S^n$ take ${\hat M}_n=D^n$, where
$D^n$ is the unit ball (disk) in $\bf R^n$ or in a $n$ dimensional
over $\bf R$ subspace in ${\cal A}_r^l$, $D_1=[0,1]$ for $n=1$. But
$S^n\setminus s_0$ has the retraction into the point in $S^n$, where
$s_0\in S^n$, $n\in \bf N$. \par Therefore, $w_0\vee
(E,f)$ and $(E,f)$ belong to the equivalence class $<(E,f)>_{t,H} :=
\{ g\in H^t_p(M,\{ s_{0,q}: q=1,...,k \} ;W,y_0): (E,f)R_{t,H}g \}
$ due to $(iii)$ and $(Y_0)$. Thus, $<w_0>_{t,H}\circ <g>_{t,H}= <g>_{t,H}$.
\par  The pseudo-manifold $M\vee
M\setminus \{ s_{0,q}\times s_{0,j}: q, j=1,...,k \} $ has the
$H^t_p$-diffeomorphism $\psi $ (see definition in \S 1.3.1) such
that $\psi (x,y)=(y,x)$ for each $(x,y)\in (M\times M\setminus \{
s_{0,q}\times s_{0,j}: q, j=1,...,k \} )$. Suppose now, that $G$ is
commutative. Then $(f\vee g)\circ \psi |_{(M\times M\setminus \{
s_{0,q}\times s_{0,j}: q, j=1,...,k \} )}= g\vee f|_{(M\times
M\setminus \{ s_{0,q}\times s_{0,j}: q, j=1,...,k \} )}$. On the
other hand, $<f\vee w_0>_{t,H}=<f>_{t,H}=<f>_{t,H}\circ
<w_0>_{t,H}=<w_0>_{t,H}\circ <f>_{t,H}$, hence, $<f\vee g>_{t,H
}=<f>_{t,H}\circ <g>_{t,H}=<f\vee w_0>_{t,H}\circ <w_0\vee g>_{t,H
}=<(f\vee w_0)\vee (w_0\vee g)>_{t,H}= <(w_0\vee g)\vee (f\vee
w_0)>_{t,H}$ due to the existence of the unit element $<w_0>_{t,H}$
and due to the properties of
$\psi $. Indeed, take a sequence $\psi _n$ as above. Therefore, the
parallel transport structure gives $(g\vee f)(\psi (x,y))=(g\circ
f)(y,x)$ for each $x, y\in M$, consequently, $(f\circ
g)R_{t,H}(g\circ f)$ for each $f, g\in H^t_p(M, \{ s_{0,q}:
q=1,...,k \} ;W,y_0)$. The using of the embedding $\chi ^*$ gives
that $(S^ME)_{t,H}$ is commutative, when $G$ is commutative.
\par The mapping $(f,g)\mapsto f\vee g$ from $H^t_p(M,\{ s_{0,q}:
q=1,...,k \} ;W,y_0)^2$  into $H^t_p(M\vee M\setminus \{
s_{0,q}\times s_{0,j}: q, j=1,...,k \}; W,y_0)$ is of class $H^t_p$.
Since the mapping $\chi ^*$ is of class $H^t_p$, then $(f,g)\mapsto
\chi ^*(f\vee g)$ is the $H^t_p$-mapping. The quotient mapping from
$H^t_p(M,\{ s_{0,q}: q=1,...,k \} ;W,y_0)$ into $(S^ME)_{t,H}$ is
continuous and induces the quotient uniformity, $T^b(S^ME)_{t,H}$
has embedding into $(S^MT^bE)_{t,H}$ for each $1\le b\le t'-t$, when
$t'>t$ is finite, for every $1\le b< \infty $ if $t'=\infty $, since
$E$ is the $H^{t'}_p$ fiber bundle, $T^bE$ is the fiber bundle with the base space $N$. Hence the multiplication
$(<f>_{t,H}, <g>_{t,H}>)\mapsto <f>_{t,H}\circ <g>_{t,H}= <f\vee
g>_{t,H}$ is continuous in $(S^ME)_{t,H}$ and is of class $H^l_p$
with $l=t'-t$ for finite $t'$ and $l=\infty $ for $t'=\infty $.

\par {\bf 4. Definition.} The $(S^ME)_{t,H}$ from Theorem
3.1 we call the wrap monoid.

\par {\bf 5. Corollary.} {\it Let $\phi : M_1\to M_2$ be a surjective
$H^t_p$-mapping of $H^t_p$-pseudo-manifolds over the same ${\cal
A}_r$ such that $\phi (s_{1,0,q})=s_{2,0,a(q)}$ for each
$q=1,...,k_1$, where $\{ s_{j,0,q}: q=1,...,k_j \} $ are marked
points in $M_j$, $j=1,2$, $1\le a\le k_2$, $l_1\le k_2$, $l_1 :=
card~ \phi (\{ s_{1,0,q}: q=1,...,k_1 \} )$. Then there exists an
induced homomorphism of monoids $\phi ^*: (S^{M_2}E)_{t,H} \to
(S^{M_1}E)_{t,H}$. If $l_1=k_2$, then $\phi ^*$ is the embedding.}
\par {\bf Proof.} Take $\Xi _1: {\hat M_1}\to M_1$ with marked points
$\{ {\hat s}_{1,0,q}: q=1,...,2k_1 \} $ as in \S 2, then take ${\hat
M}_2$ the same ${\hat M}_1$ with additional $2(k_2-l_1)$ marked
points $\{ {\hat s}_{2,0,q} : q=1,...,2k_3 \} $ such that ${\hat
s}_{1,0,q}={\hat s}_{2,0,q}$ for each $q=1,..,k_1$,
$k_3=k_1+k_2-l_1$, then $\phi \circ \Xi _1 := \Xi _2: {\hat M}_2\to
M_2$ is the desired mapping inducing the parallel transport
structure from that of $M_1$. Therefore, each ${\hat \gamma }_2:
{\hat M}_2\to N$ induces ${\hat \gamma }_1: {\hat M_1}\to N$ and to
${\bf P}_{{\hat \gamma }_2,u^2}$ there corresponds ${\bf P}_{{\hat
\gamma }_1,u^1}$ with additional conditions in extra marked points,
where $u^1\subset u^2$. The equivalence class $<(E_2,{\bf P}_{{\hat
\gamma }_2,u^2})>_{t,H}\in (S^{M_2}E)_{t,H}$ gives the corresponding
elements $<(E_1,{\bf P}_{{\hat \gamma }_1,u^1})>_{t,H}\in
(S^{M_1}E)_{t,H}$, since $DifH^t_p({\hat M}_1, \{ {\hat s}_{0,q}:
q=1,...,2k_2 \} )\subset DifH^t_p({\hat M}_1, \{ {\hat s}_{0,q}:
q=1,...,2k_3 \} )$. Then $\phi ^*(<(E_2, {\bf P}_{{\hat \gamma
}_2,u^2})\vee (E_1, {\bf P}_{{\hat \eta }_2,v^2})>_{t,H}) =\phi
^*(<(E_2, {\bf P}_{{\hat \gamma }_2,u^2})>_{t,H}) \phi ^*(<(E_1,{\bf
P}_{{\hat \eta }_2,v^2})>_{t,H})$, since $f_2\circ \phi (x)$ for
each $x\in {\Xi _1({\hat M}_1\setminus
\partial {\hat M}_1})$ coincides with $f_1(x)$, where $f_j$
corresponds to ${\bf P}_{\gamma _j,y_0\times e}$ (see also the
beginning of \S 3).
\par If $l_1=k_2$, then ${\hat M}_1={\hat M}_2$ and
the group of diffeomorphisms $DifH^t_p({\hat M}_1, \{ {\hat
s}_{0,q}: q=1,...,2k_1 \} )$ is the same for two cases, hence $\phi
^*$ is bijective and inevitably $\phi ^*$ is the embedding.

\par {\bf 6. Theorems.}  {\it {\bf 1.} There exists an alternative
topological group $(W^ME)_{t,H}$ containing the monoid
$(S^ME)_{t,H}$ and the group operation of $H^l_p$ class with
$l=t'-t$ ($l=\infty $ for $t'=\infty $). If $N$ and $G$ are
separable, then $(W^ME)_{t,H}$ is separable. If $N$ and $G$ are
complete, then $(W^ME)_{t,H}$ is complete. \par {\bf 2.} If $G$ is
associative, then $(W^ME)_{t,H}$ is associative. If $G$ is
commutative, then $(W^ME)_{t,H}$ is commutative. If $G$ is a Lie
group, then $(W^ME)_{t,H}$ is a Lie group. \par {\bf 3.} The
$(W^ME)_{t,H}$ is non-discrete, locally connected and infinite
dimensional for $dim_{\bf R}(N\times G)>1$. Moreover, if there exist
two different sets of marked points $s_{0,q,j}$ in $A_3$,
$q=1,...,k$, $j=1,2$, then two groups $(W^ME)_{t,H,j}$, defined for
$\{ s_{0,q,j}: q=1,...,k \} $ as marked points, are isomorphic.
\par {\bf 4.} The $(W^ME)_{t,H}$ has a structure of an
$H^t_p$-differentiable manifold over ${\cal A}_r$.}

\par {\bf Proof.} If $\gamma \in
H^t_p(M,\{ s_{0,q}: q=1,...,k \}; N,y_0)$, then for $u\in E_{y_0}$
there exists a unique $h_q\in G$ such that ${\bf P}_{{\hat \gamma
},u}({\hat s}_{0,q+k})=u_qh_q$, where $h_q=g_q^{-1}g_{q+k}$,
$y_0\times g_q = {\bf P}_{{\hat \gamma },u}({\hat s}_{0,q})$,
$g_q\in G$. Due to the equivariance of the parallel transport
structure $h$ depends on $\gamma $ only and we denote it by
$h^{(E,{\bf P})}(\gamma )=h(\gamma )=h$, $h=(h_1,...,h_k)$. The
element $h(\gamma )$ is called the holonomy of $\bf P$ along $\gamma
$ and $h^{(E,{\bf P})}(\gamma )$ depends only on the isomorphism
class of $(E,{\bf P})$ due to the use of $DifH^t_p({\hat M}; \{
{\hat s}_{0,q}: q=1,...,2k \} )$ and boundary conditions on $\hat
\gamma $ at ${\hat s}_{0,q}$ for $q=1,...,2k$. \par Therefore,
$h^{(E_1,{\bf P}_1)(E_2,{\bf P}_2)}(\gamma )=h^{(E_1,{\bf
P}_1)}(\gamma ) h^{(E_2,{\bf P}_2)}(\gamma )\in G^k$, where $G^k$
denotes the direct product of $k$ copies of the group $G$. Hence for
each such $\gamma $ there exists the homomorphism $h(\gamma ):
(S^ME)_{t,H}\to G^k$, which induces the homomorphism $h:
(S^ME)_{t,H}\to C^0(H^t_p(M,\{ s_{0,q}: q=1,...,k \}; N,y_0),G^k)$,
where $C^0(A,G^k)$ is the space of continuous maps from a
topological space $A$ into $G^k$ and the group structure
$(hb)(\gamma )=h(\gamma )b(\gamma )$ (see also \cite{gaj} for
$S^n$).
\par Thus, it is sufficient to construct $(W^MN)_{t,H}$
from $(S^MN)_{t,H}$. For the commutative monoid $(S^MN)_{t,H}$ with
the unit and the cancelation property there exists a commutative
group $(W^MN)_{t,H}$. Algebraically it is the quotient group $F/\sf
B$, where $F$ is the free commutative group generated by
$(S^MN)_{t,H}$, while $\sf B$ is the minimal closed subgroup in $F$
generated by all elements of the form $[f+g]-[f]-[g]$, $f$ and $g\in
(S^MN)_{t,H}$, $[f]$ denotes the element in $F$ corresponding to $f$
(see also about such abstract Grothendieck construction in
\cite{langal,swan}).

\par By the construction each point in $(S^MN)_{t,H}$ is the closed
subset, hence $(S^MN)_{t,H}$ is the topological $T_1$-space. In view
of Theorem 2.3.11 \cite{eng} the product of $T_1$-spaces is the
$T_1$-space. On the other hand, for the topological group $G$ from
the separation axiom $T_1$ it follows, that $G$ is the Tychonoff
space \cite{eng,pont}. The natural mapping $\eta : (S^MN)_{t,H }\to
(W^MN)_{t,H}$ is injective. We supply $F$ with the topology
inherited from the topology of the Tychonoff product $(S^MN)_{t,H
}^{\bf Z}$, where each element $z$ in $F$ has the form
$z=\sum_fn_{f,z}[f]$, $n_{f,z}\in \bf Z$ for each $f\in
(S^MN)_{t,H}$, $\sum_f|n_{f,z}|<\infty $. By the construction $F$
and $F/\sf B$ are $T_1$-spaces, consequently, $F/\sf B$ is the
Tychonoff space. In particular, $[nf]-n[f]\in \sf B$, hence
$(W^MN)_{t,H}$ is the complete topological group, if $N$ and $G$ are
complete, while $\eta $ is the topological embedding, since $\eta
(f+g)=\eta (f)+ \eta (g)$ for each $f, g \in (S^MN)_{t,H}$, $\eta
(e)=e$, since $(z+B)\in \eta (S^MN)_{t,H}$, when $n_{f,z}\ge 0$ for
each $f$, and inevitably in the general case $z=z^{+}-z^{-}$, where
$(z^{+}+B)$ and $(z^{-}+B)\in \eta (S^MN)_{t,H}$.
\par Using plots and $H^{t'}_p$ transition mappings of
charts of $N$ and $E(N,G,\pi ,\Psi )$ and equivalence classes
relative to $DifH^t_p(M,\{ s_{0,q}: q=1,...,k \} )$ we get, that
$(W^ME)_{t,H}$ has the structure of the $H^t_p$-differentiable
manifold, since $t'\ge t$.
\par The rest of the proof and the statements of Theorems 6(1-4)
follows from this and Theorems 3(1-3) and \cite{lujmslg,lufoclg}.
Since $(S^ME)_{t,H}$ is infinite dimensional due to Theorem 3.3,
then $(W^ME)_{t,H}$ is infinite dimensional.

\par {\bf 7. Definition.} The $(W^ME)_{t,H}
=(W^{M, \{ s_{0,q}: q=1,...,k \} }E; N,G,{\bf P})_{t,H}$ from
Theorem 6.1 we call the wrap group.

\par {\bf 8. Corollary.} {\it There exists the group homomorphism
$h: (W^ME)_{t,H}\to C^0(H^t_p(M,\{ s_{0,q}: q=1,...,k \};
N,y_0),G^k)$.}
\par {\bf Proof} follows from \S 6 and putting $h^{f^{-1}}(\gamma )=
(h^f(\gamma ))^{-1}$.

\par {\bf 9. Corollary.} {\it If $M_1$ and $M_2$ and $\phi $
satisfy conditions of Corollary 5, then there exists a homomorphism
$\phi ^*: (W^{M_2}E)_{t,H}\to (W^{M_1}E)_{t,H}$. If $l_1=k_2$, then
$\phi ^*$ is the embedding.}

\par {\bf 10. Remarks and examples.} Consider examples of $M$ which
satisfy sufficient conditions for the existence of wrap groups
$(W^ME)_{t,H}$. Take $M$, for example, $D^n_R$, $S^n_R\setminus V$
with $s_0\in
\partial V$, $D^n_R\setminus Int( D^n_b)$ with $s_0\in \partial
D^n_b$ and $0<b<R<\infty $, where $S^n_R$ denotes the sphere of the
dimension $n>1$ over $\bf R$ and radius $R$, $V$ is
$H^t_p$-diffeomorphic with the interior $Int (D^n_R)$ of the
$n$-dimensional ball $D^n_R:= \{ x\in {\bf R^n}:
\sum_{k=1}^nx_k^2\le R \} $ or in $n$ dimensional over $\bf R$
subspace in ${\cal A}_r^l$ and is the proper subset in $S^n_R:= \{
x\in {\bf R^{n+1}}: \sum_{k=1}^{n+1}x_k^2=R \}$. Instead of sphere
it is possible to take an $H^t_p$ pseudo-manifold $Q^n$ homeomorphic
with a sphere or a disk, particularly, Milnor's sphere. Indeed,
divide $M$ by the equator $ \{ x_1=0 \} $ into two parts $A_1$ and
$A_2$ and take $A_3 = \{ x\in M: x_1=0 \} \cup P$, where $s_0\in
\partial A_1\cap \partial A_2$, while $P=\emptyset $, $P=\partial
V$, $P=\partial D_b^n$ correspondingly. Then take also $V$ and
$D_b^n$ such that their equators would be generated by the equator $
\{ x_1=0 \} $ in $S_R^n$ or $D_R^n$ respectively or more generally
$Q^n$. \par Take then $M=Q^n\setminus \bigcup_{k=1}^lV_k$, where
$V_k$ are $H^t_p$-diffeomorphic to interiors of bounded quadrants
in $\bf R^n$ or in $n$ dimensional subspace in ${\cal A}_r^a$, where
$l>1$, $l\in \bf N$, $\partial V_k\cap
\partial V_j= \{ s_0 \} $ and $V_k\cap V_j=\emptyset $ for each $k\ne
j$, $diam (V_k)\le b<R/3$. In more details it is possible make a
specification such that if $l$ is even, then $[l/2]-1$ among $V_k$
are displayed above the equator and the same amount below it, two of
$V_k$ have equators, generated by equators $ \{ x_1=0 \} $ in $Q^n$.
If $l$ odd, then $[(l-1)/2]$ among $V_k$ are displayed above and the
same amount below it, one of $V_k$ has equator generated by that of
$ \{ x_1=0 \} $ in $Q^n$, $s_0\in \bigcap _k\partial V_k\cap \{ x\in
M: x_1=0 \} $. \par Divide $M$ by the equator $ \{ x_1=0 \} $ into
two parts $A_1$ and $A_2$ and let $A_3 = \{ x\in M: x_1=0 \} \cup
P$, where $P=\bigcup_{k=1}^l
\partial V_k$. Then either $A_1\setminus A_3$ and $A_2\setminus A_3$
are $H^t_p$ diffeomorphic as pseudo-manifolds or manifolds with
corners and $H^t_p$ diffeomorphic with $M\setminus [\{ s_0 \} \cup
(A_3\setminus Int (\partial A_1\cap
\partial A_2))] =: D$ or 2$(iv')$ is satisfied, since the latter
topological space $D$ is obtained from $Q^n$ by cutting a non-void
connected closed subset, $n>1$, consequently, $D$ is retractable
into a point.
\par In a case of a usual manifold $M$ the point $s_0\in \partial M$
(for $\partial M\ne \emptyset $) may be a critical point, but in the
case of a manifold with corners this $s_0$ is the corner point from
$\partial M$, since for $x\in
\partial M$ there is not less than one chart  $(U,u,Q)$ such that
$u(x)\in \partial Q$, $M\setminus \partial M=\bigcup_ku_k^{-1}(Int
(Q_k))$, $\partial M\subset \bigcup_ku_k^{-1}(\partial Q_k)$.
Further, if $M$ satisfies Conditions 2$(i-v)$ or $(i-iii,iv',v)$,
then $M\times D_R^m=P$ also satisfies them for $m\ge 1$, since
$D_R^m$ is retractable into the point, taking as two parts
$A_j(K)=A_j(M)\times D_R^m$ of $P$, where $j=1, 2$, $A_j(M)$ are
pseudo-submanifolds of $M$. Then $A_1(P)\cap A_2(P)= (A_1(M)\cap
A_2(M))\times D_R^m$ and it is possible to take $A_3(P)=A_3(M)\times
D_R^m$, $s_0(P)\in s_0(M)\times \{ x\in D_R^m: x_1=0 \} $. In
particular, for $M=S^1$ and $m=1$ this gives the filled torus.

\par This construction can be naturally generalized for non-orientable
manifolds, for example, the M\"obius band $L$, also for $M := L\setminus
(\bigcup_{j=1}^{\beta }V_j)$ with the diameter $b_j$ of $V_j$ less than the
width of $L$, where each $V_j$ is $H^t_p$ diffeomorphic with an interior of a bounded quadrant in $\bf R^2$, $s_{0,q}\in \partial L\cap (\bigcap_{j=a_1+...+a_{q-1}+1} ^{a_1+...+a_q}\partial V_j)$, $a_0 := 0$, $a_1+...+a_k=\beta $, $q=1,...,k$, since $\partial L$ is diffeomorphic with
$S^1$, also $S^1\setminus \{ s_{0,q} \}$ is retractable into a point,
consequently, $A_1$ and $A_2$ are retractable into a point. For $L$ take
${\hat M} = I^2$, then take a connected curve $\hat \eta $ consisting of
the left side $ \{ 0 \} \times [0, 1]$ joined by a straight line segment joining points $\{ 0,1 \}$ and $ \{ 1,0 \} $ and then joined by
the right side $ \{ 1 \} \times [0,1]$. This gives the proper cutting of $\hat M$ which induces the proper cutting of $L$ and of $M$ with $A_3\supset \eta \cup \partial L$ up to an $H^t_p$ diffeomorphism, where $\eta := \Xi ({\hat \eta })$, hence
the M\"obius band $L$ and $M$ satisfy Conditions 2$(i-iii,iv',v)$.

\par Take a quotient mapping $\phi : I^2\to S^1$ such that
$\phi ( \{ s_{0,1}, s_ {0,2} \} ) = s_0\in S^1$, $s_{0,1} = (0,0)$,
$s_{0,2} = (0,1)\in I^2$, where $I= [0,1]$, hence there exists the embedding $\phi ^*: (W^{S^1,s_0}E)_{t,H}\hookrightarrow (W^{I^2,
\{ s_{0,1}, s_{0,2} \} }E)_{t,H}$.

\par The Klein bottle $K$ has ${\hat M}=I^2$ with twisting equivalence relation on $\partial I^2$ so it satisfies sufficient conditions.
Moreover, $K$ is the quotient $\phi : Z\to K$ of the cylinder $Z$ with twisted equivalence relation of its ends $S^1$ using reflection relative
to a horizontal diameter. Thus $A_3\supset \phi (S^1)$. Therefore, there exists the embedding $\phi ^*: (W^{K, \{ s_0 \} }E)_{t,H} \to (W^{Z,\{ s_{0,1}, s_{0,2} \} } E)_{t,H}$, where $s_{0,1}, s_{0,2}\in \partial Z$,
$\phi ( \{ s_{0,1}, s_{0,2} \} ) = s_0$.

\par Take a pseudo-manifold $Q^n$ $H^t_p$-diffeomorphic with $S^n$
for $n\ge 2$, cut from it $\beta $ non-intersecting open domains
$V_1,...,V_{\beta }$ $H^t_p$-diffeomorphic with interiors of bounded
quadrants in $\bf R^n$,
$s_{0,q}\in \bigcap_{j=a_1+...+a_{q-1}+1} ^{a_1+...+a_q}\partial V_j$,
$a_0 := 0$, $a_1+...+a_k=\beta $, $q=1,...,k$.
Then glue for $V_1,...,V_l$, $1\le l\le \beta $, by boundaries of slits
$H^t_p$-diffeomorphic with $S^{m-1}$ the reduced product $L\vee
S^{n-2}$, since $\partial L=S^1$, $S^1\wedge S^{n-2}$ is
$H^t_p$-diffeomorphic with $S^{n-1}$ \cite{swit}. We get the
non-orientable $H^t_p$-pseudo-manifold $M$, satisfying sufficient
conditions.

\par Since the projective space ${\bf R}P^n$
is obtained from the sphere by identifying diametrically opposite
points. Then take $M$ $H^t_p$-diffeomorphic with ${\bf R}P^n$ for
$n>1$ also $M$ with cut $V_1,...,V_{\beta }$ $H^t_p$-diffeomorphic with
open subsets in ${\bf R}P^n$, $s_{0,q}\in
(\bigcap_{j=a_1+...+a_{q-1}+1} ^{a_1+...+a_q}\partial V_j)\cap \{
x\in M: x_1=0 \} $, $V_j\cap V_l=\emptyset $ for each $j\ne l$,
$a_0 := 0$, $a_1+...+a_k=\beta $, $q=1,...,k$.
Then Conditions 2$(i-v)$ or $(i-iii,iv',v)$ are also
satisfied for ${\bf R}P^n$ and $M$.
\par In view of Proposition 2.14 \cite{swit} about $H$-groups
$[X,x_0;K,k_0]$ there is not any expectation or need on rigorous
conditions on a class of acceptable $M$ for constructions of wrap
groups $(W^ME)_{t,H}$.

\par If $M_1$ is an analytic real manifold, then taking its graded
product with generators $\{ i_0,...,i_{2^r-1} \} $ of the
Cayley-Dickson algebra gives the ${\cal A}_r$ manifold (see \cite{
luoyst,lufsqv,norfamlud}). Particularly this gives $l2^r$
dimensional torus in ${\cal A}_r^l$ for the $l$ dimensional real
torus ${\bf T}_2=(S^1)^l$ as $M_1$.
\par Consider ${\bf T}_2$. It can be slit along a closed curve (loop) $C$ $H^{\infty }_p$-diffeomorphic with $S^1$ and marked points $s_{0,q}\in C\subset {\bf T}_2$ such that $C$ rotates on the surface of ${\bf
T}_2=S^1_R\times S^1_b$ on angle $\pi $ around $S^1_b$ while $C$
rotates on $2\pi $ around $S^1_R$, such that $C$ rotates on $4\pi $
around $S^1_R$ that return to the initial point on $C$, where
$0<b<R<\infty $, $q=1,...,k$, $k\in \bf N$. Therefore, the slit along $C$ of ${\bf T}_2$ is the non-orientable band which inevitably is the M\"obius band with twice larger number of marked points $ \{ s^{L}_{0,j}: j=1,...,2k \} \subset \partial L$.

\par Therefore, for $M={\bf T}_2$ as $\hat M$ take a quadrant in $\bf
R^2$ with $2k$ pairwise opposite marked points ${\hat s}_{0,q}$ and ${\hat s}_{0,q+k}$ on the boundary of $\hat M$, $q=1,...,k$, $k\in \bf N$. Suitable gluing of boundary points in $\partial \hat M$ gives the mapping $\Xi : {\hat
M}\to {\bf T}_2$, $\Xi ({\hat s}_{0,q}) = \Xi ({\hat s}_{0,q+k})= s_{0,q}$,
$q=1,...,k$. Proper cutting of $\hat M$ into ${\hat A}_j$, $j=1, 2$, or of $L$ induces that of ${\bf T}_2$. Thus we get a pseudo-submanifold $A_3({\bf T}_2) =: A_3 \supset C$, while $A_1$ and $A_2$ are retractable into a marked point $s_{0,q}\in C$ for each $q$, hence ${\bf T}_2$ satisfies Conditions 2$(i-iii,iv',v)$. In view of Corollary 9 there exists the embedding $\phi ^*: (W^{{\bf T}_2, \{ s_{0,q}: q=1,...,k \} }E)_{t,H} \to (W^{L,\{ s^{L}_{0,q}:
q=1,...,2k \} } E)_{t,H}$, where $\phi : L\to {\bf T}_2$ is
the quotient mapping with $\phi (\{ s^{L}_{0,q}, s^{L}_{0,q+k} \} ) =
\{ s_{0,q} \} $, $q=1,...,k$.

\par For the $n$-dimensional torus ${\bf T}_n$ in ${\cal A}_r^a$
with $n>2$ take a $n-1$-dimensional surface $B$ such that each its
projection into ${\bf T}_2$ is $H^t_p$-diffeomorphic with $C$ for a
loop $C$ as above. Therefore, the slit along $B$ up to a
$H^t_p$-diffeomorphism gives $M_0 := L\times I^{n-2}$ for even $n$
or $M_0 := S^1\times I^{n-1}$ for odd $n$, where $I=[0,1]$. Since
$I^m$ is retractable into a point, where $m\ge 1$. Thus we lightly get for ${\bf T}_n$ a pseudo-submanifold $A_3\supset B$ and two $A_1$ and $A_2$ retractable into points and satisfying
sufficient Conditions 2$(i-iii,iv',v)$, where ${\hat M} = I^n$ up to
a $H^t_p$-diffeomorphism, $s_{0,q}\in B\subset A_3 := A_3({\bf T}_n)$,
$\{ s^{M_0}_{0,q}, s^{M_0}_{0,q+k} \} \subset \partial M_0$, $q=1,...,k$, $k\in \bf N$. Proper cutting of $\hat M$ into ${\hat A}_j$, $j=1, 2$, induces that of ${\bf T}_n$. Thus there exists an $H^t_p$ quotient mapping $\phi : M_0\to {\bf T}_n$ with $\phi ( \{
s^{M_0}_{0,q}, s^{M_0}_{0,q+k} \} ) = \{ s_{0,q} \} $
and the embedding $\phi ^*:
(W^{{\bf T}_n, \{ s_{0,q}: q=1,...,k \} } E)_{t,H}\hookrightarrow (W^{M_0,
\{ s^{M_0}_{0,q}:  q=1,...,2k \} } E)_{t,H}$ due to
Corollary 9.

\par More generally cut from ${\bf T}_n$ open subsets $V_j$ which are $H^t_p$ diffeomorphic with interiors of bounded quadrants in $\bf R^n$ embedded into
${\cal A}_r^l$, $j=1,...,\beta $, such that $s_{0,q}\in B\cap
(\bigcap_{j=a_1+...+a_{q-1}+1} ^{a_1+...+a_q}\partial V_j)$,
$V_j\cap V_i = \emptyset $ for each $j\ne i$, $V_j\cap B = \emptyset
$ for each $j$, where $B$ is defined up to an $H^t_p$
diffeomorphism, $a_0 := 0$, $a_1+...+a_k=\beta $, $q=1,...,k$, that
gives the manifold $M_2$. Then from $M_0$ cut analogously
corresponding $V_{j,b}$, such that $s_{0,q}\in B\cap
(\bigcap_{j=a_1+...+a_{q-1}+1} ^{a_1+...+a_q}\partial V_{j,1})$,
$s_{0,q+k}\in B\cap (\bigcap_{j=a_1+...+a_{q-1}+1}
^{a_1+...+a_q}\partial V_{j,2})$, $V_{j,b_1}\cap V_{i,b_2} =
\emptyset $ for each $j\ne i$ or $b_1\ne b_2$, $a_0 := 0$,
$a_1+...+a_k=\beta $, $q=1,...,k$, $j=1,...,\beta $, $b=1, 2$, that
produces the manifold $M_1$. We choose $V_{j,b}$ such that for the
restriction $\phi : M_1\to M_2$ of the mapping $\phi $ there is the
equality $\phi (V_{j,1}\cup V_{j,2}) = V_j$ for each $j$, $\phi ( \{
s^{M_1}_{0,q}, s^{M_1}_{0,q+k} \} ) = \{ s_{0,q} \} $. This gives
the embedding $\phi ^*: (W^{M_2, \{ s_{0,q}: q=1,...,k \} }
E)_{t,H}\hookrightarrow (W^{M_1, \{ s^{M_1}_{0,q}:  q=1,...,2k \} }
E)_{t,H}$.

\par Another example is $M_3$ obtained from the previous $M_2$ with $2k$ marked points and $2\beta $ cut out domains $V_j$, when $s_{0,q}$ is identified with $s_{0,q+k}$ and each $\partial V_j$ is glued with $\partial V_{j+\beta }$ for each $j\in \lambda _q \subset \{ d: a_1+...+a_{q-1}+1\le d\le a_1+...+a_q \} $, $q=1,...,k$, $k\in \bf N$,
by an equvalence relation $\upsilon $. Such $M_3$ is obtained from
the torus ${\bf T}_{n,m}$ with $m$ holes instead of one hole in the
standard torus ${\bf T}_{n,1} ={\bf T}_n$ cutting from it $V_j$ with
$j\in \{ 1,...,2\beta \} \setminus (\bigcup_{q=1,...,k}\lambda _q)$,
where $m= m_1+...+m_k$, $m_q := card ( \lambda _q)$. For ${\bf T}_n$
and $M_2$ the surface $B$ is $H^t_p$ diffeomorphic with $(\partial
L)\times I^{n-2}$ for even $n$ or $S^1\times I^{n-1}$ for odd $n$.
Take $A_3 \supset B\cup (\bigcup_{j\in \lambda _q}\upsilon (\partial
V_j))$, it is arcwise connected and contains all marked points.
Therefore, $M_3$ satisfies conditions of \S 2 and there exists the
embedding $ \upsilon ^* : (W^{M_3, \{ s^{M_3}_{0,q}: q=1,...,k \} }
E)_{t,H}\hookrightarrow (W^{M_2, \{ s^{M_2}_{0,q}:  q=1,...,2k \} }
E)_{t,H}$. This also induces the embedding $(W^{{\bf T}_{n,m}, \{
s^{{\bf T}_{n,m} }_{0,q}: q=1,...,k \} } E)_{t,H}\hookrightarrow
(W^{{\bf T}_n, \{ s^{{\bf T}_n}_{0,q}:  q=1,...,2k-1 \} } E)_{t,H}$
such that each element $g\in (W^{{\bf T}_{n,m}, \{ s^{{\bf T}_{n,m}
}_{0,q}: q=1,...,k \} } E)_{t,H}$ can be presented as a product
$g=(..(g_1g_2)...g_m)$ of $m$ elements $g_j\in (W^{{\bf T}_n, \{
s^{{\bf T}_n}_{0,q}:  q=1,...,2k-1 \} } E)_{t,H}$, $g_j=
<f_j>_{t,H}$, $supp (\pi \circ f_j)\subset B_j$, $B_1\cup ...\cup
B_m ={\bf T}_n$, $B_i\cap B_j=\partial B_i\cap \partial B_j$ for
each $i\ne j$, each $B_j$ is a canonical closed subset in ${\bf
T}_n$, $s_{0,1}\in B_1$, $s_{0,2q}, s_{2q+1}\in B_d$ for $m_1+...+
m_{q-1} + 1 \le d \le m_1+...+m_q$, $q=1,...,k-1$, where $m_0 := 0$.

\par Evidently, in the general case for different manifolds $M$ and
$N$ wrap groups may be non isomorphic. For example, as $M_1$ take a
sphere $S^n$ of the dimension $n>1$, as $M_2$ take $M_1\setminus K$,
where $K$ is up to an $H^t_p$-diffeomorphism the union of non
intersecting interiors $B_j$ of quadrants of diameters $d_1,...,d_s$
much less, than $1$, $K=B_1\cup ... \cup B_l$, $l\in \bf N$. Let $N$
be a $\delta $-enlargement for $M_2$ in $\bf R^{n+1}$ relative to
the metric of the latter Euclidean space, where $0<\delta <\min
(d_1,...,d_l)/2$. Then the groups $(W^{M_1}N)_{t,H}$ and
$(W^{M_2}N)_{t,H}$ are not isomorphic. This lightly follows from the
consideration of the element $b:=<f>_{t,H}\in (W^{M_2}N)_{t,H}$,
where $f: M_2\to N$ is the identity embedding induced by the
structure of the $\delta $-enlargement. \par Recall, that for
orientable closed manifolds  $A$ and $B$ of the same dimension $m$
the degree of the continuous mapping $f: A\to B$ is defined as an
integer number $deg (f)\in \bf Z$ such that $f_*[A]=deg (f) [B]$,
where $[A]\in H_m(A)$ or $[B]\in H_m(B)$ denotes a generator,
defined by the orientation of $A$ or $B$ respectively
\cite{dingpan}. Consider mappings $f_j: S^n\to N$ such that
$V_j\supset
\partial B_j\cap N$, where $V_j$ is a domain in $\bf R^{n+1}$
bounded by the hyper-surface $f_j(B_j)$, $f_j$ is $w_0$ on each
$B_i$ with $i\ne j$, while the degree of the mapping $f_j$ from
$S^n$ onto $f_j(S^n)$ is equal to one. If there would be an
isomorphism $\theta : (W^{M_2}N)_{t,H}\to (W^{M_1}N)_{t,H}$, then
$\theta (b)$ would have a non trivial decomposition into the sum of
non canceling non zero additives, which is induced by mappings $f_j:
S^n\to N$. Nevertheless, an element $b$ in $(W^{M_2}N)_{t,H}$ has
not such decomposition.
\par If two groups $G_1$ and $G_2$ are not isomorphic, then
certainly $(W^ME;N,G_1,{\bf P})_{t,H}$ and $(W^ME;N,G_2,{\bf
P})_{t,H}$ are not isomorphic.

\par The construction of wrap groups can be spread on locally compact
non compact $M$ satisfying conditions 2$(ii-iv)$ or $(ii,iii,iv')$
changing $(v)$ such that $\hat M$ is locally compact non-compact
$H^t_p$-domain in ${\cal A}_r^l$, its boundary $\partial \hat M$ may
happen to be void. For this it is sufficient to restrict the family
of functions to that of with compact supports $f: M\to W$ relative
to $w_0: M\to W$, that is $supp_{w_0} (f) := cl_M \{ x\in M: f(x)\ne
y_0 \times e \} $ is compact, $cl_MA$ denotes the closure of a subset $A$ in
$M$. Then classes of equivalent elements are given with the help of
closures of orbits of the group of all $H^t_p$ diffeomorphisms $g$
with compact supports preserving marked points $DifH^t_{p,c}(M,\{
s_{0,q}: q=1,...,k \} )$ that is $supp_{id}(g) := cl_M \{ x\in M:
g(x)\ne x \} $ are compact, where $id(x)=x$ for each $x\in M$. Then
wrap groups $(W^ME)_{t,H}$ for manifolds $M$ such as hyperboloid of
one sheet, one sheet of two-sheeted hyperboloid, elliptic
hyperboloid, hyperbolic paraboloid and so on in larger dimensional
manifolds over ${\cal A}_r$. For non compact
locally compact manifolds it is possible also consider an infinite
countable discrete set of marked points or of isolated
singularities. These examples can be naturally generalized for
certain knotted manifolds arising from the given above.
\par Milnor and Lefshetz have used for $M=S^1$ and $G=\{ e \} $
the diffeomorphism group preserving an orientation and a marked
point of $S^1$. So their loop group $L(S^1,N)$ may be
non-commutative. The iterated loop group $L(S^1,L(S^{n-1},N))$ is
isomorphic with $L(S^n,N)$, where the latter group is supplied with
the uniformity from the iterated loop group, so $n$ times iterated
loop group of $S^1$ gives loop group of $S^n$ \cite{gaj}. For
$dim_{\bf R}M>1$ orientation preservation loss its significance.
Here above it was used the diffeomorphism group without any demands
on orientation preservation of $M$ such that two copies of $M$ in
the wedge product already are not distinguished in equivalence
classes and for commutative $G$ it gives a commutative wrap group.
\par Mention for comparison homotopy groups. The group $\pi _q(X)$ for
a topological space $X$ with a marked point $x_0$ in view of
Proposition 17.1 (b) \cite{botttu} is commutative for $q>1$. For
$q=1$ the fundamental group $\pi _1(X)$ may be non-commutative, but
it is always commutative in the particular case, when $X=G$ is an
arcwise connected topological group (see \S 49$(G)$ in \cite{pont}).

\par {\bf 11. Proposition.} {\it Let $L(S^1,N)$ be an $H^1_p$ loop group
in the classical sense. Then the iterated loop group
$L(S^1,L(S^1,N))$ is commutative.}
\par {\bf Proof.} Consider two elements $a, b\in L(S^1,L(S^1,N))$
and two mappings $f\in a$, $g\in b$, $(f(x))(y) = f(x,y)\in N$,
where $x, y\in I=[0,1]\subset \bf R$, $e^{2\pi x}\in S^1$. An
inverse element $d^{-1}$ of $d\in L(S^1,N)$ is defined as the
equivalence class $d^{-1}=<h^->$, where $h\in d$, $h^-(x) :=
h(1-x)$. Then \par $(1)$ $f(x,1-y)=(f(x))(1-y)\in a^{-1}$ and
$g(x,1-y)=(g(x))(1-y)\in b^{-1}$ for $L(S^1,L(S^1,N))$ and
symmetrically
\par $(2)$ $(f(y))(1-x)=f(1-x,y)\in a^{-1}$ and $(g(y))(1-x)=g(1-x,y)\in
b^{-1}$. On the other hand, $f\vee g$ corresponds to $ab$, and
$g\vee f$ corresponds to $ba$, where the reduced product $S^1\wedge
S^1$ is $H^t_p$-diffeomorphic with $S^2$ in the sense of pseudo-manifolds
up to critical subsets of codimension not less than two.

\par Consider $(S^1\vee S^1)\wedge (S^1\vee S^1)$ and $(f\vee w_0)\vee (w_0\vee g)$ and $(g\vee w_0)\vee (w_0\vee f)$ and the iterated equivalence relation $R_{1,H}$. This situation corresponds to ${\hat M}=I^2$ divided into
four quadrats by segments $\{ 1/2 \} \times [0,1]$ and $[0,1]\times \{ 1/2 \} $ with the corresponding domains for $f$, $g$ and $w_0$ in the considered wedge products, where $<f\vee w_0> = <w_0\vee f> = <f>$ is the same class of equvalent elements.

\par Since $G= \{ e \} $,
$(ab)^{-1}=b^{-1}a^{-1}$, then $g(1-x,y)\vee f(1-x,y)$ is in the
same class of equivalent elements as $g(x,1-y)\vee f(x,1-y)$. But
due to inclusions $(1,2)$ $<g(1-x,y)\vee f(1-x,y)>= <f(x,y)\vee
g(x,y)>^{-1}$ and $<f(x,y)\vee g(x,y)>= <g(x,1-y)\vee
f(x,1-y)>^{-1}$ and $<h(x,y)>^{-1}=<h(x,1-y)>=<h(1-x,y)>$ for $h\in
ab$, consequently, $<h(x,y)>=<h(1-x,1-y)>$ and $<(f\vee g)(x,1-y)> =
<f(x,1-y)\vee g(x,1-y)>\in (ab)^{-1}$, since $(x,y)\mapsto
(1-x,1-y)$ interchange two spheres in the wedge product $S^2\vee
S^2$. Hence $a^{-1}b^{-1}=b^{-1}a^{-1}$ and inevitably $ab=ba$.

\par {\bf 12. Theorem.} {\it Let $M$ and $N$ be connected both either
$C^{\infty }$ Riemann or ${\cal A}_r$ holomorphic manifolds with
corners, where $M$ is compact and $dim M\ge 1$ and $dim N>1$. Then
$(W^MN)_{t,H}$ has no any nontrivial continuous local one parameter
subgroup $g^b$ for $b\in (-\epsilon ,\epsilon )$ with $\epsilon
>0$.}
\par {\bf Proof.} Suppose the contrary, that $\{ g^b: b\in
(-\epsilon ,\epsilon ) \} $ with $\epsilon >0$ is a local nontrivial
one parameter subgroup, that is, $g^b\ne e$ for $b\ne 0$. Then to
$g^{\delta }$ for a marked $0<\delta <\epsilon $ there corresponds
$f=f_{\delta }\in H^{\infty }_p$ such that $<f>_{t,H} = g^{\delta
}$, where $f\in H^t_p$. If $f(U)= \{ y_0\times e \} $ for a sufficiently
small connected open neighborhood $U$ of $s_{0,q}$ in $M$,
then there exists a
sequence $f\circ \psi _n$ in the equivalence class $<f>_{t,H}$ with
a family of diffeomorphisms $\psi _n\in DifH^t_p(M; \{ s_{0,q}:
q=1,...,k \} )$ such that $\lim_{n\to \infty }diam \psi _n(U) =0$
and $\bigcap_{n=1}^{\infty } \psi _n(U)= \{ s_{0,q} \} $. If
$h(x)\ne y_0$, then in view of the continuity of $h$ there exists an
open neighborhood $P$ of $x$ in $M$ such that $y_0\notin h(P)$.
Consider the covariant differentiation $\nabla $ on the manifold $M$
(see \cite{kling}). The set $S_h$ of points, where $\nabla ^kh$ is
discontinuous is a submanifold of codimension not less than one,
hence of measure zero relative to the Riemann volume element in $M$.
For others points $x$ in $M$, $x\in M\setminus S_h$, all $\nabla
^kh$ are continuous.
\par Take then open $V=V(f)$ in $M$ such that
$V\supset U$ and $\nabla _{\nu }^kf|_{\partial V}\ne 0$ for some
$k\in \bf N$, where $\nabla _{\nu }f(x) := \lim_{z\to x, z\in
M\setminus V } \nabla _{\nu }f(z)$, $\nu $ is a normal (perpendicular) to $\partial V$ in $M$ at a point $x$ in the boundary $\partial V$ of $V$
in $M$.
Practically take a minimal $k=k(x)$ with such property. Since $M$ is
compact and $\partial V := cl (V)\cap cl (M\setminus V)$ is closed
in $M$, then $\partial V$ is compact. The function $x\mapsto k(x)\in
\bf N$ is continuous, since $f$ and $\nabla ^lf$ for each $l$ are
continuous. But $\bf N$ is discrete, hence each $\partial _qV:= \{
x\in
\partial V: k(x)=q \} $ is open in $V$. Therefore, $\partial V$ is a
finite union of $\partial _qV$, $1\le q\le q_m$, where $q_m :=
\max_{x\in \partial V}k(x)<\infty $ for $f=f_{\delta }$, since
$\partial V$ is compact.  Thus, there exists a subset $\lambda
\subset \{ 1,...,q_m \} $ such that $\partial V=\bigcup_{q\in
\lambda }\partial _qV$ and $\partial _qV\ne \emptyset $ for each
$q\in \lambda $. If $\nabla ^lf(x)=0$ for $l=1,...,k(x)-1$ and
$\nabla ^{k(x)}f(x)\ne 0$, then $\nabla ^{k(x)}f(\psi (y))=(\nabla
^{k(x)} f(x)).(\nabla \psi (y))^{\otimes k(x)}\ne 0$ for $y\in M$
such that $\psi (y)=x$, since $\nabla \psi (y)\ne 0$, where $\psi
\in DifH^{\infty }_p(M; \{ s_{0,q}: q=1,...,k \} )$.

\par We can take $\epsilon >0$ such that $\{ g^b: b\in (-\epsilon
,\epsilon ) \} \subset U$, where $U=-U$ is a connected symmetric
open neighborhood of $e$ in $(W^MN)_{t,H}$. Since
$g^{b_1}+g^{b_2}=g^{b_1+b_2}$ for each $b_1, b_2, b_1+b_2\in
(-\epsilon ,\epsilon )$, then $\lim_{t\to 0} g^b=e$ for the local
one parameter subgroup and in particular $\lim_{m\to \infty }
g^{1/m}=e$, where $m\in \bf N$. Take $\delta =\delta _m=1/m$ and
$f=f_m\in H^{\infty }_p$ such that $<f_m>_{t,H}=g^{1/m}$. On the
other hand, $jg^{1/m}=g^{j/m}$ for each $j<m\epsilon $, $j\in \bf
N$, hence $f_{j/m}(M) =f_{1/m}(M)$ for each $j<m\epsilon $, since
$(f\vee h) (M\vee M)=f(M)\vee h(M)$ and using embedding $\eta $ of
$(S^MN)_{t,H}$ into $(W^MN)_{t,H}$. \par The function $|\nabla
^{k(x)} _{\nu }f_{\delta }(x)|$ for $x\in
\partial V$ is continuous by $\delta $ due to the
Sobolev embedding theorem \cite{miha}, $0<\delta <\epsilon $,
consequently, $\inf_{x\in \partial V}|\nabla ^{k(x)} _{\nu
}f_{\delta }(x)|>0$, since $\partial V$ is compact. We can choose a
family $f_{\delta }$ such that $z^{(l)}(\delta ,x):=\nabla
^lf_{\delta }(x)$ is continuous for each $0\le l\le k_0$ by $(\delta
,x)\in (-\epsilon ,\epsilon )\times M$, since $\{ g^b: b\in
(-\epsilon ,\epsilon ) \} $ is the continuous by $b$ one parameter
subgroup, where $k_0:=q_m(\delta _0)$. Therefore, for this family
there exists a neighborhood $[-\epsilon +c,\epsilon -c]$ such that
$\delta _0\in [-\epsilon +c,\epsilon -c]\subset (-\epsilon ,\epsilon
)$ with $0<c<\epsilon /3$ such that $q_m(\delta )\le k_0$ for each
$\delta \in [-\epsilon +c,\epsilon -c]$ with a suitable choice of
$V(f_{\delta })$, since $\bf N$ is discrete. On the other hand,
$\sup_{x\in
\partial V(f_{\delta }), 0<\delta \le \epsilon -c} |\nabla ^{k(x)}
_{\nu }f_{\delta }(x)|\le \sup_{x\in M, 0<\delta \le \epsilon -c}
|\nabla ^{k(x)} _{\nu }f_{\delta }(x)|=: B<\infty $, since $M$ and
$[-\epsilon +c,\epsilon -c]$ are compact. \par Therefore, for this
family there exists a neighborhood $[-\epsilon +c,\epsilon -c]$ such
that $\delta _0\in [-\epsilon +c,\epsilon -c]\subset (-\epsilon
,\epsilon )$ with $0<c<\epsilon /3$ such that $q_m(\delta )\le k_0$
for each $\delta \in [-\epsilon +c,\epsilon -c]$ with a suitable
choice of $V(f_{\delta })$, since $\bf N$ is discrete. \par Then
${\underline {\lim}}_{\delta \to 0, \delta
>0} |\nabla ^{k(x)} _{\nu }f_{\delta }(x)| =: b >0$ for $x\in
\partial V$ with a suitable choice of $V=V(f_{\delta })$, since $M$
is connected, $dim M\ge 1$ and $\inf_{m\in \bf N} diam f_{j/m}(M)>0$
for a marked $\delta _0 =j/m_0<\epsilon $ with $j, m>m_0\in \bf N$
mutually prime, $(j,m)=1$, $(j,m_0)=1$. To $<f_{l/m}>_{t,H}$ there
corresponds $<f_{1/m}>_{t,H}\vee ...\vee <f_{1/m}>_{t,H
}=:<f_{1/m}>_{t,H}^{\vee l}$ which is the $l$-fold wedge product.
Thus there exists $C=const >0$ for $M$ such that $|\nabla ^{k(x)}
_{\nu }f_{l/m}(x)|\ge C l \inf_{y\in
\partial V(f_{1/m})}|\nabla ^{k(y)} _{\nu }f_{1/m}(y)|\ge Clb$,
where $C>0$ is fixed for a chosen atlas $At (M)$ with given
transition mappings $\phi _i\circ \phi _j^{-1}$ of charts. \par Consider
$\delta _0\le l/m <\epsilon -c$ and $m$ and $l$ tending to the
infinity. Then this gives $B\ge Clb$ for each $l\in \bf N$, that is
the contradictory inequality, hence $(W^MN)_{t,H}$ does not contain
any non trivial local one parameter subgroup.

\end{document}